\documentclass{article}

\usepackage{amsmath, amsthm, amsfonts}
\usepackage{color}

\usepackage{csquotes}

\newtheorem{thm}{Theorem}[section]
\newtheorem{cor}[thm]{Corollary}
\newtheorem{lem}[thm]{Lemma}
\newtheorem{theorem}[thm]{Theorem}
\newtheorem{assumption}[thm]{Assumption}

\usepackage{mdframed}
\theoremstyle{definition}
\newmdtheoremenv{boxProb}{Problem}
\newmdtheoremenv{boxDef}{Definition}
\newmdtheoremenv{boxCor}{Corollary}
\newmdtheoremenv{boxThm}{Theorem}
\newmdtheoremenv{compjob}{Computational Job}
\newmdtheoremenv{reqi}{Requirement}

\usepackage{graphicx}
\usepackage{subfigure}          

\usepackage{placeins}

\usepackage{tabularx}

\usepackage{lmodern,textcomp}

\usepackage{algorithm}
\usepackage{algpseudocode}

\usepackage{enumerate}

\usepackage{xspace} 
\newcommand\largeparbreak{\par\bigskip}

\newcommand*\tageq{\refstepcounter{equation}\tag{\theequation}}

\newcommand{\bmu}{\boldsymbol{\mu}\xspace}

\newcommand{\blambda}{\boldsymbol{\lambda}\xspace}
\newcommand{\bigeta}{\boldsymbol{\nu}\xspace}

\newcommand{\bvarrho}{\boldsymbol{\varrho}\xspace}

\newcommand{\bs}{\textbf{s}\xspace}
\renewcommand{\t}{^\textsf{T}\xspace}

\newcommand{\away}[1]{}

\newcommand{\R}{\mathbb{R}\xspace}

\newcommand{\N}{\mathbb{N}\xspace}

\newcommand{\cV}{\mathcal{V}\xspace}

\newcommand{\cB}{\mathcal{B}\xspace}

\newcommand{\cL}{\mathcal{L}\xspace}
\newcommand{\cX}{\mathcal{X}\xspace}
\newcommand{\cP}{\mathcal{P}\xspace}

\newcommand{\cJ}{\mathcal{J}\xspace}

\newcommand{\cT}{\mathcal{T}\xspace}
\newcommand{\cO}{\mathcal{O}\xspace}

\newcommand{\bB}{\textbf{B}\xspace}

\newcommand{\bD}{\textbf{D}\xspace}
\newcommand{\bS}{\textbf{S}\xspace}
\newcommand{\bP}{\textbf{P}\xspace}

\newcommand{\bx}{\textbf{x}\xspace}

\newcommand{\be}{\textbf{1}\xspace}
\newcommand{\bei}[1]{{\textbf{e}}\xspace}

\newcommand{\bI}{\textbf{I}\xspace}

\newcommand{\bO}{\textbf{0}\xspace}

\newcommand{\bu}{\textbf{u}\xspace}

\newcommand{\commentout}[1]{}

\title{High-order convergent Finite-Elements Direct Transcription Method for Constrained Optimal Control Problems}

\author{Martin P. Neuenhofen}

\begin{document}

\maketitle

\begin{abstract}
\commentout{Direct transcription with collocation polynomials on Legendre-Gauss-Radau points together with a sparse non-linear programming method is a widely used approach for numerically solving optimal control problems, which can be addressed to three reasons: First, the simplicity of implementing the collocation method. Second, the capabilities of the Radau discretization when the problem is stiff. And finally because the collocation polynomials can be chosen of arbitrary order, which yields high-order consistency. However it is known from examples that this method does not converge in the general case because it can be unstable. In particular it can diverge for bang-singular arcs and for differential-algebraic problems with a high index. Further it is undesirable that all states and controls must be discretized on the same mesh, especially when some states are much smoother or, e.g., constant.}

In this paper we present a finite element method for the direct transcription of constrained non-linear optimal control problems. 

We prove that our method converges of high order under mild assumptions. Our analysis uses a regularized penalty-barrier functional. The convergence result is obtained from local strict convexity and Lipschitz-continuity of this functional in the finite-element space.

The method is very flexible. Each component of the numerical solution can be discretized with a different mesh. General differential-algebraic constraints of arbitrary index can be treated easily with this new method. 

From the discretization results an unconstrained non-linear programming problem (NLP) with penalty- and barrier-terms. The derivatives of the NLP functions have a sparsity pattern that can be analysed and tailored in terms of the chosen finite-element bases in an easy way. We discuss how to treat the resulting NLP in a practical way with general-purpose software for constrained non-linear programming.

\end{abstract}

\section{Introduction}

\paragraph{Problem statement}
Given $n_T \in \N$, $t_0<t_1<t_2<...<t_{(n_T-1)}=t_E \in \R$ and $\Omega := (t_0,t_E)$. We consider the following optimal control problem (OCP).
\begin{subequations}
\begin{align}
	& \operatornamewithlimits{min}_{(y,z) \in \cX} 
					& F(y,z) = \int_\Omega f\big(\,\dot{y}(t),y(t),z(t),t\,\big)&\ \mathrm{d}t						                     \label{eqn:OCP1}\\
	& \text{s.t.} 	& b\big(\,y(t_0),y(t_1),...,y(t_E)\,\big)&= 0\,,     \label{eqn:OCP2}\\
	& 				& c\big(\,\dot{y}(t),y(t),z(t),t\,\big)  &= 0 
	\quad \text{f.a.e. } t \in \Omega\,,             \label{eqn:OCP3}\\
	& 				& z(t)                                   &\geq 0
	\quad \text{f.a.e. } t \in \Omega                        \label{eqn:OCP4}
\end{align}\label{eqn:OCP}
\end{subequations}
In this problem statement there appears the following Hilbert space
\begin{align*}
	\cX &:= \left(H^1(\Omega)\right)^{n_y} \times \left(L^2(\Omega)\right)^{n_z}\,,\\
	\langle(y,z),(v,w)\rangle_\cX &:= \sum_{j=1}^{n_y} \langle y_{[j]},v_{[j]} \rangle_{H^1(\Omega)} + \sum_{j=1}^{n_z} \langle y_{[j]},v_{[j]} \rangle_{H^1(\Omega)}
\end{align*}
for $n_y,n_z \in \N_0$, $n_x := n_y+n_z$, and there are the following functions
\begin{align*}
	f &: \R^{n_y} \times \R^{n_y} \times \R^{n_z} \times \Omega \rightarrow \R\\
	b &: \underbrace{\R^{n_y} \times ... \times \R^{n_y}}_{\text{$n_T$ times}} \rightarrow \R^{p}\\
	c &: \R^{n_y} \times \R^{n_y} \times \R^{n_z} \times \Omega \rightarrow \R^{m}
\end{align*}
for $m,p \in \N_0$.

$F$ in \eqref{eqn:OCP1} is a cost-functional that has to be minimized by the function $x=(y,z)$ subject to the constraints \eqref{eqn:OCP2}--\eqref{eqn:OCP4}. \eqref{eqn:OCP2} defines punctual conditions, e.g. boundary conditions, whereas \eqref{eqn:OCP3} are called differential-algebraic equations (DAE) \cite{kunkel2006differential}. $y$ are differential functions, since their derivatives appear in the objective and the constraints. $z$ are called auxiliary variables. They are time-dependent but may be non-smooth or may be discontinuous. Finally, \eqref{eqn:OCP4} gives inequality constraints. \enquote{f.a.e.} means \enquote{for almost every} in the sense of Lebesgue-measurability. This is because in contrast to elements of the Sobolev space $H^1(\Omega)$ the value in one time-point of a function in $L^2(\Omega)$ is meaningless due to lack of smoothness \cite{adams2003sobolev}.

In this paper we provide a direct transcription scheme for \eqref{eqn:OCP} and prove its convergence.

\paragraph{Numerical solution of optimal control problems}

Numerical methods for solving optimal control problems work by discretizing the infinite-dimensional solution space into a vector space of finite dimension. They can be classified into two times two types: Direct and indirect transcription methods, and among each a full and a semi-discretization approach. The first distinction says whether either the equations of the problem itself are discretized (direct) or the equations of its optimality system --- known as Karush-Kuhn-Tucker (KKT) equations \cite{Gander2014} --- are discretized (indirect). The semi-discrete approach does not fully discretize the solution space $\cX$ but determines several components of $x$ (called \textit{states}) from others (called \textit{controls}), where only the controls are discretized and the states are determined from integrating the DAE for the numerically determined controls. The DAE can be integrated with standard numerical methods \cite{Hairer:1993:SOD:153158}. The semi-discrete approach only works in rare cases because a separation into dependent and independent variables of $x$ is not possible in general.

Direct methods are preferred over indirect methods in practice . This is for a simple reason: The indirect approach leads to an equation system, whose solutions are potential minimizers. Such general equation systems can only be solved numerically for a very accurate initial guess (some variant of Newton's method \cite{Nash2001265}), which is usually unavailable. The direct approach instead leads to a non-linear program, whose local minimizers can be always obtained from feasible initial guesses (e.g. with a feasible-path SQP method, cf. \cite{Jian:2005:FDS:1716511.1716843}), which are usually available.

A survey on optimal control problems and numerical methods to solve them and references on all the aforementioned can be found in \cite{rao_asurvey,Betts1998,Conway2012}.

\paragraph{Practical methods} GPOPS-II \cite{Patterson:2014:GMS:2684421.2558904} is a widely used software to solve (all problems that can be expressed as) \eqref{eqn:OCP} in a reformulated way. It is a direct full discretization method that discretizes $\cX$ by piecewise polynomial functions. The discretization of the equations is done as follows: \eqref{eqn:OCP3} is discretized by collocation, i.e. it is only satisfied in a finite number of values $t\in\Omega$. The integral \eqref{eqn:OCP1} is solved by quadrature on the collocation points with the Gauss-Legendre weights on the Radau points. The collocation method used is called Radau 2A \cite{hairer2010solving} or Legendre-Gauss-Radau \cite{Wang:2012:LCM:2317881.2317893}.

Direct transcription with collocation polynomials on Legendre-Gauss-Radau points together with a sparse non-linear programming method is a widely used approach for numerically solving optimal control problems, which can be addressed to three reasons: First, the simplicity of implementing the collocation method. Second, the capabilities of the Radau discretization when the problem is stiff. And finally because the collocation polynomials can be chosen of arbitrary order, which yields high-order consistency. However, it is known from examples that the direct transcription with Gauss-Legendre-Radau collocation does not necessarily converge for all optimal control problems. In particular it is divergent for problems with singular optimal controls, cf. \cite[p.1566]{KAMESWARAN20061560} and \cite{ALYCHAN1973}. To the best of our knowledge, it remains an open problem to prove that this method converges for problems that have general path constraints, this is, problems of a form as considered in \eqref{eqn:OCP}. Further it is undesirable --- since potentially leading to unnecessarily large NLPs --- that all states and controls must be discretized on the same mesh, especially when some states are much smoother that others.

\paragraph{Scope of the paper} In contrast to the aforementioned short-comings of other methods, in this paper we present a direct transcription method that has the following advantages: First, for our method we prove high-order convergence under mild assumptions. Second, different finite-element spaces can be used for each component $x_{[j]}$, $j=1,...,n_x$. This allows to use highly tailored search spaces for each component of the solution function. 

If, for example, $x_{[1]}(t)$ models a temperature and $x_{[2]}$ the vertical position of a wheel on a rough ground, then our method could use a grid for $x_{[1]}(t)$ that is much coarser than for $x_{[2]}$. This capability can also be used to allow only particular shapes for some components (e.g. piecewise constant or piecewise linear functions). 

Finally, our method treats the equality constraints in a way that makes it irrelevant of which type they are. To clarify this: Here we present the method with point-constraints and DAE constraints, but at the end of the day these constraints are combined into a residual that is simply treated by a penalty term. As long as this residual is well-defined for the elements in the solution space $\cX$ it does not matter of which kind the constraints are. 

\paragraph{Structure} In Section~2 we introduce regularity assumptions that we make on the problem. Then we introduce some parameters for our analysis.
In Section~3 we replace \eqref{eqn:OCP} by an unconstrained variational problem with penalty- and barrier-terms.
In Section~4 we introduce the finite element method that we use to compute a local minimizer of the unconstrained variational problem.
Section~5 describes how the method can be implemented and how the resulting non-linear program can be solved with general-purpose NLP solvers.

\section{Parameters of the problem and the method}

We make some assumptions to ensure that the problem is well-posed. Then we introduce five parameters that are needed to define the method. Afterwards we define theoretical parameters that will be needed only for the convergence analysis.

\paragraph{Assumptions on the problem (problem-parameters)}
We define the squared constraint residual norm $r$ for elements $x=(y,z)\in\cX$:
\begin{align*}
r(y,z) := \|c\big(\dot{y}(\cdot),y(\cdot),z(\cdot),\cdot\big)\|^2_{L^2(\Omega)} + \|b\big(y(t_0),...,y(t_E)\big)\|_2^2
\end{align*}
We sometimes omit the word \enquote{squared} in the following.

Say that for problem \eqref{eqn:OCP} we are interested in the particular local minimizer $x^\star \in \cX$. We assume the existence of parameters $0<\delta \in \R$, $L_f,L_r \in \R$, $\pi \in \R$, $F_{\text{min}}\in\R$, $F_{\text{opt}} \in \R$, such that the following holds for $\cB := \lbrace u \in \cX \, \vert \, \|x^\star-u\|_\cX \leq \delta \rbrace$:
\begin{subequations}
\begin{align}
	F(x) &\geq F(x^\star) & & \forall x \in \cB\\
	F_{\min} &\leq \operatornamewithlimits{\min}_{\tilde{x}\in\cX}\lbrace\,F(\tilde{x})\,\rbrace\\
	F_{\text{opt}} & \geq F(x^\star)\\
	\|x^\star\|_\cV &\leq \pi\\
	|F(u)-F(v)| &\leq L_f \, \|u-v\|_\cV & &\forall u,v \in \cB\\
	|r(u)-r(v)| &\leq L_r^2\, \|u-v\|_\cV^2 & & \forall u,v \in \cB
\end{align}
\end{subequations}
The first condition says that $x^\star$ is a local minimizer in a neighbourhood of radius $\delta$. The second and third relations make requirements on the boundedness of the cost-functional. The fourth condition bounds of the minimizer itself. And the latter two Lipschitz-conditions make sure that the sensitivities of optimality gap and feasibility residual are bounded locally around $x^\star$.

\paragraph{Method-parameters}
At this stage we assume that the user provides three parameters: A mesh-size $0<h \in \R$, a mesh-ratio coefficient $0<\sigma \in \R$, and a degree $d \in \lbrace 0,1,2,...,30\rbrace$ that defines the polynomial degree of the finite-element functions. From these user-defined numerical values we define a penalty- and a barrier-parameter.
\begin{align*}
	\omega := h^{d/2}\,,\quad\quad\tau := h^{d}
\end{align*}

\paragraph{Theoretical parameters for the analysis}
The following definitions of parameters help our analysis at a later stage, but their actual values do not need to be known in practice. We use the Bachmann-Landau notation \cite{Knuth:1976:BOB:1008328.1008329}. Assuming that the problem-parameters do all live in $\Theta(1)$ and $0<h \in \cO(1)$, we provide for each parameter below how it scales with $h,\tau,\omega\rightarrow +0$.
\begin{subequations}
\begin{align}
	\delta_\omega 	&:= \delta_\pi & & \in \Theta(1)\\
	L_\omega 		&:= L_f + \frac{\omega\,\delta_\omega}{2} + \frac{1}{2\,\omega}\,L_r^2\,\delta_\pi & & \in \Theta\left(\omega+\frac{1}{\omega}\right)\\
	\pi_\omega 		&:= \frac{2}{\omega}\,\sqrt{F_{\min}+F_{\text{opt}}} & & \in \Theta\left(\frac{1}{\omega}\right)\\
	\alpha_\omega 	&:= \frac{\omega}{2} & & \in \Theta(\omega)\\[7pt]
	\delta_{\omega,\tau} 	&:= \min\left\lbrace\,\delta_\pi\,,\,\frac{\sigma\,h}{d+1}\,\frac{\tau}{2\,L_\omega}\,\right\rbrace & & \in \Theta\left(h\,\tau\,\left(\omega+\frac{1}{\omega}\right)\right)\\
	\eta &:= \frac{\tau}{2\,L_\omega} & &\in \Theta\left(\tau\,\left(\omega+\frac{1}{\omega}\right)\right)\\
	L_{\omega,\tau} 		&:= L_\omega + \frac{t_E-t_0}{\eta}\,\tau\,n_z\,\frac{d+1}{\sigma\,h} & & \in \Theta\left( \left(1+\frac{1}{h}\right)\,\left(\omega+\frac{1}{\omega} \right) \right)  \\
	\pi_{\omega,\tau} 		&:= \pi_{\omega} + \sqrt{\tau\,n_x\,\frac{t_E-t_0}{\alpha_\omega}} & & \in\Theta\left( \frac{1}{\omega}+ \sqrt{\frac{\tau}{\omega}\,} \right) \\
	\alpha_{\omega,\tau} 	&:= \frac{\omega}{2} & & \in \Theta\left( \omega \right)
\end{align}
\end{subequations}

\section{Unconstrained minimization approach}
In this section we approximate the constrained original problem \eqref{eqn:OCP} by an unconstrained problem. We first introduce a replacement for the equality- and then for the inequality-constraints. We provide an analysis on how the local minimizers of the respective problems are related to each other.

\paragraph{Penalty form}
From \eqref{eqn:OCP} we derive the following related problem for a parameter $\omega>0$.
\begin{subequations}
	\begin{align*}
	& \operatornamewithlimits{min}_{(y,z) \in \cX} 
	& F_\omega(y,z) := 	& F(y,z) + \frac{\omega}{2}\,\|(y,z)\|_\cX^2 + \frac{1}{2\,\omega}\, r(y,z) \tageq\label{eqn:rOCP1}\\[7pt]
	&\text{s.t.} 	& 	z(t)                                   \geq& 0
	\quad \text{f.a.e. } t \in \Omega                        \tageq\label{eqn:rOCP2}
	\end{align*}\label{eqn:rOCP}
\end{subequations}
Due to the convexification by a square of $\|(y,z)\|_\cX$ it follows that for each local minimum of $F$ the functional $F(\cdot)+\frac{\omega}{2}\,\|\cdot\|_\cX^2$ has a strict minimum, required $\omega$ is sufficiently small.

In general one cannot prove convergence for a minimizer $x_\omega^\star$ of \eqref{eqn:rOCP} towards a minimizer $x^\star$ of \eqref{eqn:OCP}. But, one can show that the characterizing properties of a minimizer, namely the optimality gap and the constraint residual norm, converge to zero. The following theorem shows this.

\begin{theorem}[Convergence of the penalty solution]\label{Thm:Convergence_rOCP}
	Given a strict local minimizer $x^\star_\omega$ of \eqref{eqn:rOCP}. Then one of the local minimizers $x^\star$ of \eqref{eqn:OCP} satisfies:
	\begin{align*}
	F(x^\star_\omega)-F(x^\star) &\leq \frac{\omega}{2}\,\|x^\star\|_\cX^2 \in \cO(h^{d/2}) \\
	r(x^\star_\omega) &\leq 2\,\omega\,\big(F_{\text{opt}}-F_{\min}\big)+\omega^2\,\|x^\star\|_\cX^2 \in \cO(h^{d/2}) \\
	\|x_\omega^\star\|_\cX &\leq \sqrt{\frac{2}{\omega}\,\big(F_{\text{opt}}-F_{\min}\big)\,} \in \cO(h^{-d/4})
	\end{align*}
	I.e., there is convergence of the optimality gap and the feasibility residual. Further, $\|x_\omega^\star\|_\cX$ is bounded.
\end{theorem}
\begin{itshape}
	\noindent\textbf{Proof:} Since $x_\omega^\star$ is optimal for \eqref{eqn:rOCP} and $z_\omega^\star,z^\star\geq 0$ are both feasible for \eqref{eqn:rOCP}, it holds:
	\begin{align}
		F_\omega(x^\star_\omega) \leq F_\omega(x^\star)\label{eqn:aux:ConvReg}
	\end{align}
	
	For the first proposition, we use $F(x^\star_\omega) \leq F_\omega(x^\star_\omega)$ and $F_\omega(x^\star) = F(x^\star) + \frac{\|x^\star\|_\cX^2}{2}\,\omega$. Inserting both into \eqref{eqn:aux:ConvReg} yields the result.
	
	For the second proposition we multiply \eqref{eqn:aux:ConvReg} by $2\,\omega$.
	\begin{align*}
	r(x^\star_\omega) + 2\,\omega\,F(x^\star_\omega) + \omega^2\,\|x^\star_\omega\|_\cX^2 \leq 2\,\omega\,F(x^\star)+ \omega^2\,\|x^\star\|_\cX^2
	\end{align*}
	Since non-negative, we omit $\omega^2\,\|x^\star_\omega\|_\cX^2$ on the left-hand side.
	Then, moving all terms except the constraint residual norm to the right, we find:
	\begin{align*}
	r(x^\star_\omega)\leq 2\,\omega\,\big(F(x^\star)-F(x^\star_\omega)\big) + \omega^2\,\|x^\star\|_\cX^2
	\end{align*}
	
	For the third proposition, we use $F_{\text{min}}+\frac{\omega}{2}\,\|x_\omega^\star\|^2_\cX \leq F_\omega(x_\omega^\star)$. Inserting this into \eqref{eqn:aux:ConvReg} yields
	\begin{align*}
		\frac{\omega}{2}\,\|x_\omega^\star\|^2_\cX \leq -F_{\text{min}} + F(x^\star)\,.
	\end{align*}
	q.e.d.
\end{itshape}

At this stage we have defined an analytical strict local minimizer $x^\star_\omega$ for a regularized problem \eqref{eqn:rOCP}. We also have investigated the sense in which $x^\star_\omega$ is related to a local minimizer $x^\star$ of the original problem \eqref{eqn:OCP}.

\paragraph{Penalty-barrier form}
In \eqref{eqn:rOCP} there are inequality constraints. It is a common approach from interior-point methods \cite{Potra:2000:IM:365036.365063,Gondzio_interiorpoint} to use barrier functions in order to replace a constrained problem by an unconstrained one. Using logarithmic barriers with barrier-parameter $\tau>0$ for the components of $z = (z_{[1]},z_{[2]},...,z_{[n_z]})$, we reformulate \eqref{eqn:rOCP} into the following:
\begin{align}
	& \operatornamewithlimits{min}_{(y,z) \in \cX} 
	& F_{\omega,\tau}(y,z) := 	& F_\omega(y,z) - \tau\,\sum_{j=1}^{n_y} \int_\Omega \log\big(z_{[j]}(t)\big)\ \mathrm{d}t\label{eqn:brOCP}
\end{align}
The local minimizers of \eqref{eqn:brOCP} we call $x^\star_{\omega,\tau}$.
To analyze how this problem is related to \eqref{eqn:rOCP}, we first introduce a general result on the barrier-approach. 

\begin{theorem}[Convergence of barrier solutions]\label{Theorem:BarrierConvergence}
Let $\Omega = [t_0,t_E] \subset \R$, $\cV$ a Hilbert space, $\cV \subset \left(L^2(\Omega)\right)^n$, $n \in \N$. Given $F : \cV \rightarrow \R$. Consider the problem
\begin{align}
	\min_{u \in \cV} F(u)\text{ , s.t. }u_{[j]}(t)\geq 0\text{ f.a.e. }t \in \Omega,\ j \in \cJ\subset\lbrace 1,...,n\rbrace\label{eqn:thm:BarrierConvergence:gen_noBarrier}
\end{align}
with a strict local minimizer $u^\star \in \cV$.
Choose $0<c\in \R$, $0 < \delta \in \R$, $0<\tau \in \R$, and define $\cB:= \lbrace v \in \cV \, \vert \, \|u^\star - v \|_\cV\leq \delta\rbrace$. Consider the problem 
\begin{align}
\min_{u \in \cV} F_\tau(u):=F(u) - \tau\,\sum_{j\in\cJ} \int_\Omega \log\big(u_{[j]}(t)\big)\ \mathrm{d}t\label{eqn:thm:BarrierConvergence:gen_Barrier}
\end{align}
with a strict local minimizer $u_\tau^\star \in \cV$.

If the following conditions are satisfied
\begin{enumerate}[(i)]
	\item $F$ is convex in $\cB$
	\item $F(v) \geq F(u^\star) + c \,\|u^\star-v\|_\cV^2\quad \forall v \in \cB$
	\item $u_\tau^\star \in \cB$
\end{enumerate}
then it holds:
\begin{align}
	\|u_\tau^\star - u^\star\|_\cV \leq \sqrt{\tau\,n\,\frac{t_E-t_0}{c}\,}
\end{align}
\end{theorem}
\begin{itshape}
	\noindent\textbf{Proof:} 
	Since $F$ is convex in $\cB$, we can use:
	\begin{align*}
		F(u^\star) \geq F(u^\star_\tau) - \delta F(u^\star_\tau)(u_\tau^\star-u^\star)
	\end{align*}
	From optimality we know
	\begin{align*}
		\delta F_\tau(u^\star_\tau)(v) \equiv \delta F(u^\star_\tau)(v) - \tau\,\sum_{j\in\cJ} \int_\Omega \frac{v_{[j]}(t)}{u^\star_{\tau,[j]}(t)}\ \mathrm{d}t=0\quad \forall v \in \cV
	\end{align*}
	Inserting this into the above for $v = u_\tau^\star-u^\star$ yields:
	\begin{align}
		F(u^\star) \geq F(u^\star_\tau) - \tau\,\sum_{j=1}^{n} \int_\Omega \frac{u^\star_{\tau,[j]}(t)-u^\star_{[j]}(t)}{u^\star_{\tau,[j]}(t)}\ \mathrm{d}t\label{eqn:aux:int_gap}
	\end{align}
	We use that
	\begin{align*}
		0 \leq \frac{u^\star_{\tau,[j]}(t)-u^\star_{[j]}(t)}{u^\star_{\tau,[j]}(t)} \leq 1\quad \text{f.a.e. }t \in \Omega\,,
	\end{align*}
	which in turn is equivalent to $u^\star_{\tau,[j]}(t) \geq u^\star_{[j]}(t)$. This is obviously true because $\delta F(u^\star)(\cdot)=0$ and the gradient of the barrier pushes $u_\tau^\star$ in an amount $\geq 0$ into positive direction relative to $u^\star$, i.e. $u^\star_\tau\geq u^\star$.
	
	We insert the bound into \eqref{eqn:aux:int_gap} and obtain:
	\begin{align*}
	F(u^\star_\tau) \leq F(u^\star) + \tau\,\sum_{j\in\cJ} \int_\Omega \underbrace{\frac{u^\star_{\tau,[j]}(t)-u^\star_{[j]}(t)}{u^\star_{\tau,[j]}(t)}}_{\leq 1}\ \mathrm{d}t
	\end{align*}
	From requirement, we can replace the left-hand side with its lower bound $F(u^\star) + c \,\|u^\star_\tau-u^\star\|_\cV^2$. Insertion yields
	\begin{align*}
		F(u^\star) + c \,\|u^\star_\tau-u^\star\|_\cV^2 \leq F(u^\star) + \tau\,|\cJ|\,|\Omega|
	\end{align*}
	q.e.d.
\end{itshape}

We give some remarks on condition $(iii)$. One can enforce it by enlarging $\delta$ through the addition of a strictly convex term to $F$, or alternatively by choosing $\tau$ sufficiently small. In general, a condition like $(iii)$ is unavoidable because too large values of $\tau$ would always push $u_\tau^\star$ out of $\cB$.
\largeparbreak

We apply Theorem~\ref{Theorem:BarrierConvergence} with \eqref{eqn:rOCP} for \eqref{eqn:thm:BarrierConvergence:gen_noBarrier} and \eqref{eqn:brOCP} for \eqref{eqn:thm:BarrierConvergence:gen_Barrier}. The following corollary states that the requirements of the theorem are satisfied.
\begin{cor}
	Choose $h>0$ sufficiently small. Define $\cB := \lbrace u \in \cX \, \vert \, \|u - x^\star_\omega\|_\cX \leq \delta_\omega \rbrace$. Then the following holds:
	\begin{enumerate}[(i)]
		\item $F_\omega$ is convex in $\cB$.
		\item $F_\omega(v) \geq F(x^\star) + c_\omega \, \|x^\star - v\|^2_\cV \quad \quad \forall v \in \cB$
		\item $x^\star_{\omega,\tau} \in \cB$
	\end{enumerate}
\end{cor}
\begin{itshape}
	\noindent\textbf{Proof:} $(i)$ follows from the definition of $\delta_\omega$ and the requirements on the convexity of $F$, which implies local convexity of $F_\omega$. $(ii)$ follows from the definition of $F_\omega$ and $c_\omega$. Finally, $(iii)$ can be shown in the following way: For the limit of $\tau\rightarrow+0$ we know that $F_{\omega,\tau}$ is strictly convex in $\cB$ and has a strict local minimizer in $\cB$. For $\tau>0$ we know that convexity is still maintained. There is thus a value $\tau>0$ sufficiently small such that the strict minimizer of $F_{\omega,\tau}$ remains in $\cB$. As $\tau\rightarrow 0$ for $h \rightarrow 0$, there is a value $h>0$ sufficiently small such that condition $(iii)$ is satisfied. q.e.d.
\end{itshape}

Using the theorem we find:
\begin{align}
	\|x^\star_{\omega,\tau}-x^\star_\omega\|_{\cX} \leq \sqrt{\tau\,n_z\,\frac{t_E-t_0}{c_\omega}\,} \in \cO(h^{d/4})
\end{align}

\paragraph{Positivity of $z_{\omega,\tau}^\star$}
In this paragraph we show that the component $z_{\omega,\tau}^\star$ of $x^\star_{\omega,\tau}$ is strictly positive with a positive lower bound. To this end we start from the following general result.
\begin{theorem}[Strict positivity of barrier solutions]\label{thm:strictPos}
Consider the problem \eqref{eqn:thm:BarrierConvergence:gen_Barrier} with the strict local minimizer $u_\tau^\star$ for $0<\tau \in \R$. Let $0<L\in\R$ be bounded, such that $\|\delta F(u_\tau^\star)(\cdot)\|_{\cV'}\leq L$ is satisfied. Then the following holds:
\begin{align*}
	u^\star_{\tau,[j]}(t)\geq \frac{\tau}{L} \quad \text{ f.a.e. }t\in\Omega\quad \forall j \in \cJ\,.
\end{align*}
\end{theorem}
\begin{itshape}
	\noindent\textbf{Sketch of the proof:} We consider a one-dimensional example. Consider for $u : \Omega=(t_0,t_E) \rightarrow \R$ the functional
	\begin{align*}
		F_\tau(u) := \underbrace{\int_\Omega L\,u(t) \ \mathrm{d}t}_{=:F(u)}\ -\ \tau\,\int_\Omega \log\big( u(t) \big) \ \mathrm{d}t
	\end{align*}
	The larger the constant $L$ is, the more is $u$ pushed to the boundary. The Euler-Lagrange-ODE yields the unique strict minimizer $u^\star_\tau(t)=\frac{\tau}{L}$, which is a constant function. $L$ is the upper bound of the slope of $F$. The above example is the worst-case scenario where the slope is oriented so to maximally force $u(t)$ towards zero. We conclude that from the above example follows the worst case in terms of a lower bound for the pointwise value of $u_\tau^\star(t)$.
\end{itshape}

Applied to \eqref{eqn:brOCP} the theorem says the following:
\begin{cor}
	Consider $x_{\omega,\tau}^\star=(y_{\omega,\tau}^\star,z_{\omega,\tau}^\star)$. Then it holds:
	\begin{align*}
	z^\star_{\omega,\tau,[j]}(t)\geq \frac{\tau}{L_\omega}\quad\text{ f.a.e. }t \in \Omega\quad\quad\forall j=1,...,n_z\,.
	\end{align*}
\end{cor}
\begin{itshape}
	\noindent\textbf{Proof:} follows from Theorem~\ref{thm:strictPos}.
\end{itshape}

\section{Finite Element Method (FEM)}

We want to compute a numerical approximation $x^\star_h$ to $x_{\omega,\tau}^\star$. We start with introducing a finite-element space in which we search $x_h^\star$. Then we propose a general convergence result for finite-element methods for non-linear variational problems. Finally, we apply the convergence result to our numerical approximation. We will find $\|x^\star_h-x_{\omega,\tau}^\star\|_\cX \in \cO(h^{(d-3)/2})$.


\paragraph{Finite element space and discrete variational problem}

In this section we describe a finite-dimensional Hilbert space $\cX_h \subset \cX$. 
The following technicalities define the construction of the function spaces.

For an interval $T \subset \Omega$ we write $\cP_d(T)$ for the space of polynomials $p : T \rightarrow \R$ of degree $d \in \N_0$. We call $\hat{T}:=(0,1)$ (not necessarily a subset of $\Omega$) the unit interval. We define triangulations. $\cT_h = \lbrace T_j \subset \Omega,\ j=1,...,n_{\cT_h}\rbrace$ is called triangulation if $T_j$ are disjunct intervals whose union is $\Omega$. We write $|T_j|$ for the length of the interval $T_j$. We use the mesh-size $h>0$ of a mesh, that satisfies $h \geq |T_j|$ $\forall j=1,...,n_{\cT_h}$. We require quasi-uniformity, i.e. $\min_{j,k}\frac{|T_j|}{|T_k|} \geq \sigma$ for a constant $\sigma \in \R$.

Given $n_x$ triangulations $\cT_h^{(1)},\cT_h^{(2)}...,\cT_h^{(n_x)}$ and a degree $d \in \N_0$, we define the space
\begin{align*}
	\cX_{h,d} := \Big\lbrace x=(y,z) \in \cX \ \Big\vert x_{[j]} \in \cP_d(T_k)\ \forall\, T_k\in \cT_h^{(j)} \ \forall j=1,...,n_x \Big\rbrace\,.
\end{align*}
Notice that since $x \in \cX$, it follows that $y$ is continuous, whereas $z$ is not necessarily continuous. Notice further that for each component of $x$ we can use a different triangulation $\cT_h$. As discussed in the introduction, this can be of crucial benefit in practice since smoother functions can be discretized on coarser meshes.

Using the above finite element space $\cX_{h}$ (where the sub-index of $d$ is omitted sometimes) we discretize \eqref{eqn:brOCP} into the following unconstrained non-linear programming problem:
\begin{align}
& \operatornamewithlimits{min}_{(y,z) \in \cX_h} 
& F_{\omega,\tau}(y,z) := 	& F_\omega(y,z) - \tau\,\sum_{j=1}^{n_y} \int_\Omega \log\big(z_{[j]}(t)\big)\ \mathrm{d}t\label{eqn:dbrOCP}
\end{align}
We denote the local minimizers of \eqref{eqn:dbrOCP} with $x_h^\star$.

\largeparbreak

In the following we discuss the approximation error in $\cX_{h,d}$. From the Bramble-Hilbert-Lemma we find
\begin{align}
	&\operatornamewithlimits{inf}_{v \in \cX_{h,d}}\lbrace\,\|x_{\omega,\tau}^\star - v\|_\cX \,\rbrace \\
	\leq & n_x \, C_d\, \sigma\, h^{d-1} \, (t_E-t_0) \, \max_{1\leq k \leq n_x}\Bigg\lbrace\ \max_{T_j \in \cT_h^{(k)}} \Big\lbrace \, \frac{\|x_{\omega,\tau,[k]}^\star\|_{H^{d}(T_j)}}{|T_j|} \,\Big\rbrace  \ \Bigg\rbrace\,,\label{eqn:BrambleHilbert}
\end{align}
where $C_d$ is a constant that depends only on the degree $d$. The bound is not immediately useful because on the right-hand side there appear Sobolev-$(d,p)$-norms of the components of $x_{\omega,\tau}^\star$ over the intervals of the discretization. In order for these norms to remain bounded for large degrees $d$, we need to make the assumption that $x_{\omega,\tau}^\star$ is infinitely smooth on each interval $T_j$. We formulate this assumption more precisely.

\begin{assumption}[Smoothness of $x_{\omega,\tau}$]
	We make the following assumptions for every $T_j \in \bigcup_{1\leq l \leq n_x}\cT_h^{(k)}$ for every $k \in \lbrace 1,...,n_x\rbrace$ and every $d \in \N_0$:
	\begin{enumerate}[(I)]
		\item \textbf{Smoothness of the original solution} We make the following assumption on the problem instance: $\|x_{[k]}^\star\|_{H^{d}(T_j)} \in \cO(|T_j|)$. This means that the original problem instance has a solution $x^\star$ for which triangulations $\cT_h^{(k)}$ can be found such that in each interval $x^\star$ is infinitely smooth.
		\item \textbf{Smoothness of the penalty-barrier solution} We assume $\|x_{\omega,\tau,[k]}^\star\|_{H^{d}(T_j)} \in \cO(\|x^\star_{[k]}\|_{W^{d,p}(T_j)})$. We consider this assumption reasonable because the penalty with $\omega$ smooths the derivatives and the equality constraints, and the barrier with $\tau$ smooths the inequality constraints. This is why in general we expect $x_{\omega,\tau}^\star$ to be smoother than $x^\star$.
	\end{enumerate}
\end{assumption}

Using the two above assumptions we can bound the right-hand side in \eqref{eqn:BrambleHilbert} and find:
\begin{align}
	\operatornamewithlimits{inf}_{v \in \cX_{h,d}}\lbrace\,\|x_{\omega,\tau}^\star - v\|_\cX \,\rbrace \in \cO(h^{d-1}) \label{eqn:ConvergenceInf}
\end{align}

\paragraph{Lipschitz-continuity of $F_{\omega,\tau}$ in $\cX_{h,d}$}
The functional $F_{\omega,\tau}$ in \eqref{eqn:dbrOCP} uses logarithmic barriers to force strict positivity of the $z$-components of the minimizer. To apply a convergence result, we need Lipschitz-continuity of $F_{\omega,\tau}$ in a spherical neighbourhood around $x^\star_{\omega,\tau}$. 

The paragraph is organized as follows. We start with a norm equivalence of the infinity-norm and $2$-norm for elements of $\cX_{h,d}$. Using this equivalence, we then show the Lipschitz-continuity.
\begin{thm}[Equivalence of norms]\label{thm:EquivalenceOfNorms}
Consider $\cX_{h,d}$ for $d \in \N_0$, $d\leq 30$, $h > 0$. Then the following holds:
\begin{align*}
	\max_{1\leq k \leq n_x}\,\|u_{[k]}\|_{L^\infty(\Omega)} \leq \frac{d+1}{\sigma\,h}\,\|u\|_\cX\quad  \forall u \in \cX_{h,d}
\end{align*}
\end{thm}
\begin{itshape}
	\noindent\textbf{Proof:}
	Consider a component $v := u_{[k]}$ of an arbitrary function $u \in \cX_{h,d}$. On each interval of $T_j \in \cT_h^{(k)}$ the function $v$ is a polynomial of degree $d$. 
	
	Using a computer algebra software, we solved the following convex quadratic optimization problem in $a_1,...,a_d \in \R$:
	\begin{subequations}
		\begin{align}
			& 								& \hat{v}(t) &:= 1 + \sum_{j=1}^d a_j \, t^j\\
			&\min_{a_1,...,a_j \in \R}		& G(\hat{v}) &:=\|\hat{v}\|_{L^2(\hat{T})}^2
		\end{align}
	\end{subequations}
	We observed that for a respective value $d \in \lbrace 0,1,2,...,30\rbrace$ the following holds for the minimizer $\hat{v}(t)$:
	\begin{align*}
		\|\hat{v}\|_{L^\infty(\Omega)} &= \hat{v}(0) = 1\\
		\|\hat{v}\|_{L^2(\Omega)} & = \frac{1}{d+1}
	\end{align*}
	Notice that $\hat{v}$ has been chosen such that the value $\|\hat{v}\|_{L^2(\Omega)}$ is minimized. Since norms are linear to the scaling of $\hat{v}$, we can conclude the following general result:
	\begin{align}
		\|\hat{v}\|_{L^2(\hat{T})} \geq \frac{\|\hat{v}\|_{L^\infty(\hat{T})}}{d+1} \quad \forall \hat{v} \in \cP_d(\hat{T}) \label{eqn:aux:EquivalenceNorm_NumericalResult}
	\end{align}
	
	Since $\Omega$ is the union of all $T_j \in \cT_h^{(k)}$, it follows:
	\begin{align*}
		\exists j\in\lbrace 1,...,n_{\cT^{(k)}_h}\rbrace\quad : \quad \|v\|_{L^\infty(T_j)} = \|v\|_{L^\infty(\Omega)}
	\end{align*}
	For this particular $j$ we can use \eqref{eqn:aux:EquivalenceNorm_NumericalResult}, i.e.:
	\begin{align*}
		\|v\|_{L^2(T_j)} \geq \frac{\|v\|_{L^\infty(T_j)}}{d+1}\,|T_j|
	\end{align*}
	Using $|T_j| \geq \sigma \,h$ and $\|v\|_{L^\infty(T_j)}=\|v\|_{L^\infty(\Omega)}$ yields the proposition.
	q.e.d.
\end{itshape}

We now show the Lipschitz-continuity of $F_{\omega,\tau}$.
\begin{lem}[Penalty-barrier Lipschitz constant]\label{lem:BarrierLipschitz}
	Define $\cB:= \lbrace (y,z) \in \cX \, \vert \, \|(y,z) - x_{\omega,\tau}^\star\|_\cX \leq \delta_{\omega,\tau} \rbrace$. Then it holds:
	\begin{align}
	|F_{\omega,\tau}(y,z)-F_{\omega,\tau}(v,w)| \leq L_{\omega,\tau}\,\|(y,z)-(v,w)\|_\cX \quad \forall (y,z),(v,w) \in \cB\cap\cX_h
	\end{align}
\end{lem}
\begin{itshape}
	\noindent\textbf{Proof:}
	\begin{align*}
	&|F_{\omega,\tau}(y,z)-F_{\omega,\tau}(v,w)| \\
	\leq &\underbrace{|F_{\omega,}(y,z)-F_{\omega,\tau}(v,w)|}_{\leq L_\omega \, \|(y,z)-(v,w)\|_\cX} + \tau \,\sum_{j=1}^{n_z}\int_\Omega |\log\big(z_{[j]}(t)\big)-\log\big(w_{[j]}(t)\big)|\ \mathrm{d}t
	\end{align*}
	Since the first term is known to be L-continuous, we now focus on the second term. For some component $j$ it holds:
	\begin{align}
		&\tau \,\sum_{\tilde{j}=1}^{n_z}\int_\Omega |\log\big(z_{[\tilde{j}]}(t)\big)-\log\big(w_{[\tilde{j}]}(t)\big)|\ \mathrm{d}t \\
		\leq &\tau\,n_z\,\int_\Omega |\log\big(z_{[j]}(t)\big)-\log\big(w_{[j]}(t)\big)|\ \mathrm{d}t
		\label{eqn:Lem:PenBarLipschitz}
	\end{align}
	We know that $z^\star_{\omega,\tau,[j]}(t)\geq \frac{\tau}{L_\omega}$ and $(y,z),(v,w) \in \cB\cap\cX_{h,d}$. We apply Lemma~\ref{thm:EquivalenceOfNorms}, from which we obtain the following bounds on the gap from $z^\star_{\omega,\tau,[j]}$ to $z_{[j]}$ and $w_{[j]}$:
	\begin{align*}
		\|z^\star_{\omega,\tau,[j]} - z_{[j]}\|_{L^\infty(\Omega)} \leq \frac{d+1}{\sigma\,h}\,\delta_{\omega,\tau}
	\end{align*}
	Using the bound (also for $w_{[j]}$) we find
	\begin{align*}
		\min\lbrace\,z_{[j]}(t)\,,\,w_{[j]}(t)\,\rbrace \geq \frac{\tau}{L_{\omega}} - \frac{d+1}{\sigma\,h}\,\delta_{\omega,\tau} \equiv \eta\quad\text{ f.a.e. }t \in \Omega
	\end{align*}
	From Theorem~\ref{thm:EquivalenceOfNorms} we have:
	\begin{align*}
		\|z_{[j]} - w_{[j]}\|_{L^\infty(\Omega)} \leq \frac{d+1}{\sigma\,h}\,\|z_{[j]} - w_{[j]}\|_{L^2(\Omega)}
	\end{align*}
	We can bound the integrand as
	\begin{align*}
				&|\log\big(z_{[j]}(t)\big)-\log\big(w_{[j]}(t)\big)|  \\
		\leq 	&|\frac{\mathrm{d}}{\mathrm{d}s}\log(s)\Big\vert_{s=\eta}|\,\|z_{[j]} - w_{[j]}\|_{L^\infty(\Omega)}\\
		\leq 	& \frac{1}{\eta} \, \frac{d+1}{\sigma\,h}\,\underbrace{\|z_{[j]} - w_{[j]}\|_{L^2(\Omega)}}_{\leq \|(y,z)-(v,w)\|_\cX}\,.
	\end{align*}
	Inserting the above bound for the integrand into the right-hand side of \eqref{eqn:Lem:PenBarLipschitz} yields the proposition. q.e.d.
\end{itshape}

\paragraph{General convergence result}

The following theorem provides a general convergence result for discretized unconstrained variational problems.
\begin{thm}\label{Thm:ConvergenceFEM}
	Let $\cV$ be a Hilbert space, $F : \cV \rightarrow \R$ a functional, $u^\star \in \cV$ a local minimizer of $F$. Let $0 < \delta \in \R$, $0 < c$, $C < \infty$ $c,C \in \R$, and define
	\begin{align*}
		\cB := \lbrace v \in \cV \ \vert \ \|u^\star - v\|_\cV \leq \delta \rbrace\,.
	\end{align*}
	Let $\cV_h \subset \cV$ be a Hilbert space. Define
	\begin{align*}
		u_h^\star := \operatornamewithlimits{argmin}_{v \in \cV_h\cap\cB}\lbrace\,F(v)\,\rbrace\,.
	\end{align*}
	If the following conditions hold
	\begin{enumerate}[(i)]
		\item $F(u^\star) + c\,\|u^\star - v\|_\cV^2 \leq F(v) \leq F(u^\star) + C\,\|v - u^\star\|_\cV$\quad $\forall v \in \cB \cap \cV_h$
		\item $\operatornamewithlimits{inf}_{v \in \cV_h}\lbrace\|u^\star-v\|_\cV\rbrace \leq \delta$
	\end{enumerate}
	then $u^\star_h$ exists, is unique, and satisfies:
	\begin{align*}
		\|u^\star - u_h^\star\|_\cV \leq \sqrt{\frac{C}{c}}\ \sqrt{\operatornamewithlimits{inf}_{v \in \cV_h}\lbrace\|u^\star-v\|_\cV\rbrace}
	\end{align*}
\end{thm}
\begin{itshape}
	\noindent\textbf{Proof:} By requirement $(ii)$, the best-approximation $\tilde{u}_h \in \cV_h$ to $u^\star$ lives in $\cB$. We can insert it into $(i)$.
	\begin{align*}
		F(u^\star) + C\,\|\tilde{u}_h - u^\star\|_\cV \geq F(\tilde{u}_h) \geq F(u^\star_h) \geq F(u^\star) + c\,\|u^\star_h - u^\star\|_\cV^2
	\end{align*}
	q.e.d.
\end{itshape}

The theorem makes two requirements: $(i)$ needs strict convexity of the cost-functional in the neighbourhood of the local minimizer with a constant $c>0$. Also, Lipschitz-continuity with a constant $C$ in the local region $\cB$ is required. The second condition $(ii)$ enforces that the finite-element space holds an element that is close to the minimizer and lives in the local convex region $\cB$.
\largeparbreak

In the following we apply Theorem~\ref{Thm:ConvergenceFEM} to \eqref{eqn:dbrOCP} and \eqref{eqn:brOCP} to show convergence of $x^\star_h$ towards $x^\star_{\omega,\tau}$. The following lemma gives the result.
\begin{lem}
	Let $h>0$ be sufficiently small and $d \in \lbrace 4,...,30\rbrace$. Then:
	\begin{align*}
		\|x_{\omega,\tau}^\star - x_h^\star\|_{\cX} \in \cO(h^{(d-3)/2})
	\end{align*}
\end{lem}
\begin{itshape}
	\noindent\textbf{Proof:} 
	$F_{\omega,\tau}$ satisfies the condition $(i)$ of Theorem~\ref{Thm:ConvergenceFEM} with $\delta=\delta_{\omega,\tau}$, $c=\alpha_{\omega,\tau}$ and $C=L_{\omega,\tau}$. Also condition $(ii)$ is satisfied for $h>0$ sufficiently small because $\delta_{\omega,\tau} \in \Theta(h^{(2+d)/2})$ whereas $\operatornamewithlimits{inf}_{v \in \cV_h}\lbrace\|u^\star-v\|_\cV\rbrace \in \cO(h^{d-1})$.
	
	We find
	\begin{align*}
		\|x_{\omega,\tau}^\star - x_h^\star\|_{\cX} \leq \sqrt{\frac{L_{\omega,\tau}}{\alpha_{\omega,\tau}}}\,\operatornamewithlimits{inf}_{v \in \cV_h}\lbrace\|u^\star-v\|_\cV\rbrace \in \cO(h^{(d-3)/2})\,.
	\end{align*}
	q.e.d.
\end{itshape}

Finally, we combine all the convergence results that we had so far.
\begin{thm}[Convergence]
Use $d \in \lbrace 4,...,30\rbrace$ and $h>0$ sufficiently small. Then every local minimizer of $x_h^\star$ of \eqref{eqn:NLP} there is a local minimizer $x^\star$ of \eqref{eqn:OCP} such that the following hold:
\begin{align*}
	F(x^\star_h) - F(x^\star) &\in \cO(h^{(d-3)/4})\\
	r(x_h^\star) & \in \cO(h^{(d-3)/4})
\end{align*}
\end{thm}
\begin{itshape}
	\noindent\textbf{Proof:} Using the triangular inequality with $\|x_h^\star-x^\star_{\omega,\tau}\|_\cX \in \cO(h^{(d-3)/2})$ and $\|x_\omega^\star-x^\star_{\omega,\tau}\|_\cX \in \cO(h^{d/4})$, we find
	$\|x_h^\star-x^\star_{\omega}\|_\cX \in \cO(h^{(d-3)/4})$. From the Lipschitz-continuity of $F$ and $r$ follows:
	\begin{align*}
		|F(x^\star_h) - F(x^\star_\omega)| &\leq L_f \,\|x_h^\star-x^\star_{\omega,\tau}\|_\cX \in \cO(h^{(d-3)/4})\\
		|r(x^\star_h) - r(x^\star_\omega)| &\leq L_r^2 \,\|x_h^\star-x^\star_{\omega,\tau}\|^2_\cX \in \cO(h^{d-3}/2) \subset \cO(h^{(d-3)/4})
	\end{align*}
	q.e.d.
\end{itshape}

Under \enquote{sufficient growth}-conditions on $r$ and $F$ one could prove that $\|x_\omega^\star-x^\star\|_\cX$ converges in $\cO(\omega)$. In particular, one would need sufficient growth of $F$/ $r$ when a point $x$ leaves the feasible/ optimal region, respectively. Using another triangular inequality, this would imply:
\begin{align*}
	\|x_h^\star-x^\star\|_\cX \in \cO(h^{(d-3)/4})
\end{align*}
Remark: When $d\leq 3$ then one could still prove convergence by choosing the size of $\omega,\tau$ differently with regard to $h$. This is orthogonal to the scope of this paper because we are particularly interested in large values of $d$, i.e. high-order convergence.

\section{Implementation}
In this section we discuss how the non-linear programming problem \eqref{eqn:dbrOCP} can be implemented efficiently, such that a parametrization vector $\bx \in \R^n$ for $x_h^\star$ can be computed with a non-linear programming software, e.g., WORHP \cite{Bueskens2013}, SNOPT \cite{Gill:2005:SSA:1055334.1055404}, IPOPT \cite{Waechter2006}, or Knitro \cite{ByrdNoceWalt06}. Though, we find it crucial to point out that in the general  case non-linear programming problems do not admit a time-efficient solution (in terms of polynomial complexity) \cite{Pardalos1991}, except $P=NP$ holds \cite{Fortnow:2009:SPV:1562164.1562186}.

To describe the implementation we assume that triangulations $\cT_h^{(k)}$ are given. The end-points of all intervals $T_j$ from the triangulations $\cT_h^{(k)}$ of all components $k=1,...,n_x$ defines a new triangulation $\cT_{\text{quad}}$, which consists of intervals between all neighbouring end-points. We define a quadrature rule by using Gauss-Legendre polynomials of degree $2\,d$ on each interval of $\cT_{\text{quad}}$. In result we obtain quadrature weights $\alpha_j$ and quadrature points $\rho_j \in \Omega$ for $j=1,...,M$, such that an integral over an arbitrary sufficiently smooth function $g$ can be approximated as
\begin{align*}
	\int_\Omega g(t)\,\mathrm{d}t = \sum_{j=1}^M \alpha_j\,g(\rho_j) + \cO(h^{2\,d})\,.
\end{align*}

\paragraph{Notation}

We define the following vectors for arbitrary $\ell \in \N$:
\begin{align*}
\vec{\alpha}_\ell &:= (\underbrace{\alpha_1,...,\alpha_1}_{\ell\text{ times}},\underbrace{\alpha_2,...,\alpha_2}_{\ell\text{ times}},...,\underbrace{\alpha_M,...,\alpha_M}_{\ell\text{ times}})\t \in \R^{(\ell\,M)}
\end{align*}

We introduce the notation of block-diagonal matrices. Given matrices $\bB_\ell \in \R^{J \times K}$ for $\ell=1,...,L$. We define
\begin{align*}
\operatornamewithlimits{blockdiag}_{\ell=1,...,L}(\bB) \in \R^{(J\,L) \times (K\,L)}
\end{align*}
as the block-diagonal matrix of $L$ diagonal blocks, where the $\ell$th diagonal block is $\bB_\ell$.

To evaluate functions like $f\big(\dot{y}(t),y(t),z(t),t\big)$ on the point $t=\rho_j$, we use the short writing $f_j$.

The given triangulations $\cT_h^{(k)}$ define $\cX_{h,d}$, where $\dim(\cX_{h,d})=:N \in \N$. We can define matrices
\begin{align*}
	\bP_{\dot{y},y,z} 	&\in \R^{\left((2\,n_y+n_z)\,M\right) \times N}\\
	\bP_{y(t_\star)}  	&\in \R^{(n_y\,n_{T}) \times N}\\
	\bP_\varrho 			&\in \R^{M \times 1}\\
	\bP_\lambda 		&\in \R^{(m\,M)\times (m\,M)}\\
	\bP_\nu 			&\in \R^{p \times p}\\
	\bP_\mu 			&\in \R^{(n_z\,M)\times(n_z\,M)}\,.
\end{align*}
We explain the mapping of these matrices: For an element $x \in \cX_{h,d}$ we can use a parametric representation $\bx \in \R^N$. We define $\bP_{\dot{y},y,z}$ such that it maps $\bx$ to $(\dot{y},y,z)$ evaluated in the quadrature points $\rho_j$.
\begin{align}
	\begin{pmatrix}
	\dot{y}(\rho_1)\\
	y(\rho_1)\\
	z(\rho_1)\\
	\dot{y}(\rho_2)\\
	y(\rho_2)\\
	z(\rho_2)\\
	\vdots\\
	\dot{y}(\rho_M)\\
	y(\rho_M)\\
	z(\rho_M)	
	\end{pmatrix} = \bP_{\dot{y},y,z}\, \bx\label{eqn:evaluation_vector}
\end{align}
Analogously, we define $\bP_{y(t_\star)}$ as the pointwise evaluations of $y$ at the points $t_0,t_1,...,t_E$.
\begin{align*}
\begin{pmatrix}
y(t_0)\\
y(t_1)\\
\vdots\\
y(t_{E})
\end{pmatrix} = \bP_{y(t_\star)}\, \bx
\end{align*}
For ease of presentation we define the other matrices in the following way:
\begin{align*}
	\bP_\varrho 		&:= \vec{\alpha}_1\\
	\bP_\lambda 	&:= \operatornamewithlimits{blockdiag}_{j=1,...,M}(\sqrt{\alpha_j}\,\bI_{m \times m})\\
	\bP_\nu 		&:= \bI_{p \times p}\\
	\bP_\mu 		&:= \bI_{(n_z\,M)\times(n_z\,M)}
\end{align*}
We define a matrix $\bS \in \R^{N \times N}$ as 
\begin{align*}
	\bS := \bP_{\dot{y},y,z}\t\,\operatornamewithlimits{blockdiag}_{j=1,...,M}(\bI_{(n_y+n_x)}\,\alpha_j)\,\bP_{\dot{y},y,z}\,.
\end{align*}
For a function $x=(y,z)\in\cX$ we define $g_j := \big(z(\rho_j)\big)^{\alpha_j}$ for $j=1,...,M$, where the power of the vector $z(\rho_j) \in \R^{n_z}$ is meant for each component.


\paragraph{Rewriting the non-linear program}
We can rewrite \eqref{eqn:dbrOCP} as follows,
\begin{align}
&\min_{\bx \in \R^N} & & F(\bx) + \frac{\omega}{2}\,\|\bx\|_\bS^2 + \frac{1}{2\,\omega}\,\|H(\bx)\|^2_2 + \tau\,\|\log\big(G(\bx)\big)\|_1\label{eqn:NLP}
\end{align}
where we used the following functions,
\begin{align*}
	F(\bx) 		&:= \sum_{j=1}^M \alpha_j\,f_j\\
	H_c(\bx) 	&:= \begin{pmatrix}
						\sqrt{\alpha_1}\,c_1\\
						\sqrt{\alpha_2}\,c_2\\
						\vdots\\
						\sqrt{\alpha_M}\,c_M
					\end{pmatrix}\\
	H_b(\bx) 	&:= b(\bP_{y(t_\star)}\,\bx)\\
	H(\bx) 		&:= \begin{pmatrix}
						H_c(\bx)\\
						H_b(\bx)
					\end{pmatrix}\\
	G(\bx) 		&:= \begin{pmatrix}
						g_1\\
						g_2\\
						\vdots\\
						g_M
					\end{pmatrix}\,,
\end{align*}
where $f_j,c_j,g_j$ are meant in $\dot{y}(\rho_j),{y}(\rho_j),{z}(\rho_j),\rho_j$ as computed from $\bx$.

We can rewrite some of these functions:
\begin{align*}
F(\bx) 		&= \bP_\varrho\t	\,\operatornamewithlimits{blockdiag}_{j=1,...,M}(f_j)\,\be\\
H_c(\bx) 	&= \bP_\lambda\t\,\operatornamewithlimits{blockdiag}_{j=1,...,M}(c_j)\,\be\\
G(\bx) 		&= \bP_\mu\t 	\,\operatornamewithlimits{blockdiag}_{j=1,...,M}(g_j)\,\be
\end{align*}
We further define the Lagrangian function:
\begin{align*}
	\cL : \, &\R^N \times \R^1 \times \R^{(m\,M)}\times \R^p \times \R^{n_z\,M} \rightarrow \R\,,\ \\ &(\bx;\bvarrho,\blambda,\bigeta,\bmu) \mapsto\\ &\bvarrho\t\,\Big(F(\bx)+\frac{\omega}{2}\,\|\bx\|_\bS^2\Big)-\blambda\t\,H_c(\bx)-\bigeta\t\,H_b(\bx)-\bmu\t\,G(\bx)
\end{align*}

\paragraph{Derivatives}
We write $JF,JH,JH,JG$ for the Jacobians of $F,H,G$ with respect to $\bx$, and $\nabla^2_\bx\cL(\bx;\bvarrho,\blambda,\bigeta,\bmu)$ for the Hessian of $\cL$ with respect to $\bx$. We derive formulas for these derivatives.

As a notational trick, we interpret 
\begin{align*}
	\begin{pmatrix}
		\varrho(\rho_1)\\
		\varrho(\rho_2)\\
		\vdots\\
		\varrho(\rho_M)
	\end{pmatrix} = \bP_\varrho\,\bvarrho,\quad 
	\begin{pmatrix}
		\lambda(\rho_1)\\
		\lambda(\rho_2)\\
		\vdots\\
		\lambda(\rho_M)
	\end{pmatrix} = \bP_\lambda\,\blambda,\quad
	\begin{pmatrix}
		\mu(\lambda_1)\\
		\mu(\lambda_2)\\
		\vdots\\
		\mu(\lambda_M)
	\end{pmatrix} = \bP_\mu\,\bmu,
\end{align*}
where $\varrho : \Omega \rightarrow \R$, $\lambda : \Omega \rightarrow \R^m$, $\mu : \Omega \rightarrow \R^{n_z}$. The trick simplifies the presentation of the derivatives. It is orthogonal to our scope to analyse these "functions" any further.

We use $Jf_j$, $Jc_j$ and $Jg_j$ for the Jacobians of $f,c,g$ at $\rho_j$ with respect to $(\dot{y},y,z) \in \R^{n_y} \times \R^{n_y} \times \R^{n_z}$. We find:
\begin{align*}
JF(\bx) 	&= \bP_\varrho\t	\,\operatornamewithlimits{blockdiag}_{j=1,...,M}(Jf_j)\,\bP_{\dot{y},y,z}\\
JH_c(\bx) 	&= \bP_\lambda\t\,\operatornamewithlimits{blockdiag}_{j=1,...,M}(Jc_j)\,\bP_{\dot{y},y,z}\\
JH_b(\bx) 	&= \bP_\eta\t 	\,Jb(\bP_{y(t_\star)}\,\bx)\,\bP_{y(t_\star)}\\
JG(\bx) 	&= \bP_\mu\t 	\,\operatornamewithlimits{blockdiag}_{j=1,...,M}(Jg_j)\,\bP_{\dot{y},y,z}
\end{align*}
We define $\ell_j := \varrho(\rho_j)\,f_j-\lambda(\rho_j)\t\,c_j-\mu(\rho_j)\t\,g_j$ for $j=1,...,M$, where $\varrho,\lambda,\mu$ are computed from $\bvarrho,\blambda,\bmu$. Using the notation $\nabla^2\ell_j$ for the Hessian of $\ell_j$ with respect to $(\dot{y},y,z)$, we can write:
\begin{align*}
	\nabla^2_\bx\cL(\bx;\bvarrho,\blambda,\bigeta,\bmu) = \,&\bvarrho\,\omega\,\bS \\
	&+\bP_{\dot{y},y,z}\t\,\operatornamewithlimits{blockdiag}_{j=1,...,M}(\nabla^2\ell_j)\,\bP_{\dot{y},y,z}
	\\ &-\bP_{y(t_\star)}\t\,\left[\nabla_{[y(t_0),...,y(t_E)]}^2\big(\bigeta\t\,b(\bP_{y(t_\star)}\,\bx)\big)\right]\,\bP_{y(t_\star)}
\end{align*}

We discuss the sparsity pattern of the derivative matrices. They consist of products of the form
\begin{align*}
	\bP_a\t \, \bD \, \bP_b
\end{align*}
where $\bP_a,\bP_b$ can be interpreted as finite-element matrices that map parametrization vectors into evaluation-vectors like \eqref{eqn:evaluation_vector} and where $\bD$ is a block-diagonal matrix, whose $j$th block is the evaluation of a derivative of $f,c,\ell$ at the $j$th quadrature point. Depending on the finite-element bases, which are reflected by the columns of $\bP_a$ and $\bP_b$, the product inherits the structure of $\bD$. If $\bP_a,\bP_b$ are banded, then so is the product $\bP_a\t \, \bD \, \bP_b$.

\paragraph{Solution of the NLP with general-purpose solvers}
Different approaches are possible for solving \eqref{eqn:NLP}. Most trivial, one can apply a numerical method for unconstrained minimization to compute a local minimizer $\bx$. This approach appears undesirable to pursue because the gradient of the cost-function may change rapidly with a perturbation in $\bx$, which is due to the terms with $H(\bx)$ and $G(\bx)$. This can result in that the method selects very small steps per iteration and thus requires a lot of iterations in total, which is computationally prohibitive.

There are primal interior-point methods that are more tailored to \eqref{eqn:NLP}. To quote from \cite{ForsgrenGill}: \enquote{Standard second-order line-search or trust-region methods are easily adapted to find a local minimizer of [problem \eqref{eqn:NLP}] (for details see, e.g., Dennis and Schnabel \cite{dennis+schnabel:1983}).} Actually the method from \cite{ForsgrenGill} and references therein are designed for solving constrained non-linear programs of the form
\begin{subequations}
\begin{align}
	&\operatornamewithlimits{min}_{\bx \in \R^N}& F(\bx)& \\
	&\text{s.t.}								& H(\bx)&=\bO\\
	& 											& G(\bx)&\geq\bO\,.
\end{align}\label{eqn:constrNLP}
\end{subequations}
But they treat them by minimizing a regularized primal merit-function that has precisely the structure of \eqref{eqn:NLP}. A minimizer $\bx$ to this merit-function is found by solving \eqref{eqn:NLP} iteratively from initial guesses that are computed as minimizers of \eqref{eqn:NLP} for larger parameters of $\omega,\tau$. For further literature on interior-point methods we refer to \cite{Potra:2000:IM:365036.365063,DBLP:books/daglib/0000928,Nocedal06} and the references therein.

\largeparbreak
In the following we consider the situation where a general-purpose numerical software for the solution of constrained non-linear programs shall be applied to solve \eqref{eqn:NLP}. To raise chances that the software finds a minimizer in a small amount of iterations, it is beneficial to avoid strongly non-linear algebraic expressions both in the cost-function and the constraints.

It seems desirable to compute $\bx$ by solving a problem like \eqref{eqn:constrNLP}, because $H$ and $G$ appear only linearly and the scales of $F,H$ and $G$ are separated from each other. It is well-known that for $\omega,\tau\rightarrow 0$ the solution of \eqref{eqn:NLP} converges to that of \eqref{eqn:constrNLP}, which is why obtaining $\bx$ from \eqref{eqn:constrNLP} appears so attractive at first glance \cite{Waechter2006,ForsgrenGill}. Unfortunately, this is a wrong conclusion as we demonstrate: While the solution of \eqref{eqn:NLP} converges to that of \eqref{eqn:constrNLP}, the same does not hold the other way round. This is because for our problem the functions $F,H,G$ do actually depend on $\omega,\tau$, as these depend on $h$. Most easily the implications can be seen when using a quadrature rule that yields a large number of points $M$. When $M$ exceeds $N$ then probably the problem \eqref{eqn:constrNLP} becomes infeasible, whereas \eqref{eqn:NLP} obviously always remains feasible. The only thing that would change in \eqref{eqn:NLP} is that $\|H_c(\bx)\|_2^2$ would become a more accurate approximation of $r(y,z)$.

Instead of using \eqref{eqn:constrNLP} we propose to solve a problem that is more equivalent to \eqref{eqn:NLP}, but that uses the format of a constrained non-linear programming problem. The formulation we suggest is given below.
\begin{align*}
	&\operatornamewithlimits{min}_{\bx,{\footnotesize \begin{pmatrix}
		\blambda\\
		\bigeta
		\end{pmatrix}},\bs}	& 		
		& F(\bx)+\frac{\omega}{2}\,\|\bx\|_\bS^2 + \frac{\omega}{2}\,\left\|\begin{pmatrix}
		\blambda\\
		\bigeta
		\end{pmatrix}\right\|_2^2\\
	& \text{s.t.} 												&
		& H(\bx) - \omega\,\begin{pmatrix}
		\blambda\\
		\bigeta
		\end{pmatrix} = \bO\\
	& & & G(\bx) - \bs = \bO\\
	& & & \bs \geq \bO
\end{align*}
In this problem the penalties with $\omega$ have been reformulated equivalently, which yields equality constraints that are linear in $H$ and $\blambda,\bigeta$. The large penalty in the cost-function has been replaced by a small convex quadratic term. The barriers with $\tau$ have been reformulated in a way that is no more equivalent to \eqref{eqn:NLP}. That is because it is not possible to replace logarithmic barriers by a combination of linear and quadratic constraints. However, when this problem is passed to an interior-point method with logarithmic barriers, such as IPOPT or Knitro, then this NLP method will actually solve for some small $\vartheta>0$ the problem
\begin{subequations}
\begin{align}
&\operatornamewithlimits{min}_{\bx,{\footnotesize \begin{pmatrix}
		\blambda\\
		\bigeta
		\end{pmatrix}},\bs}	& 		
& F(\bx)+\frac{\omega}{2}\,\|\bx\|_\bS^2 + \frac{\omega}{2}\,\left\|\begin{pmatrix}
\blambda\\
\bigeta
\end{pmatrix}\right\|^2 + \vartheta\,\|\log(\bs)\|_1\\
& \text{s.t.} 												&
& H(\bx) - \omega\,\begin{pmatrix}
\blambda\\
\bigeta
\end{pmatrix} = \bO\,,\\
& & & G(\bx) - \bs = \bO\,.
\end{align}\label{eqn:IPOPT_eqn}
\end{subequations}
If by an algorithmic option it is possible to set $\vartheta=\tau$ then finally the above problem is again equivalent to \eqref{eqn:NLP}. E.g., in IPOPT this option can be set as follows.
\begin{align*}
	\texttt{options.mu\_min := }&\tau\\
	\texttt{options.mu\_target := }&\tau
\end{align*}
In this case, IPOPT will compute a sequence of solutions to problem \eqref{eqn:IPOPT_eqn} for a decreasing sequence of values for $\vartheta>0$, until finally a solution is computed and returned for the value $\vartheta=\tau$.

\section{Conclusion}
We presented a finite element method for the solution of optimal control problems with general DAE constraints, pointwise constraints and inequality constraints. We provided a rigorous proof for the convergence of the method toward feasible locally optimal points.

Two questions are left open. First, how does the method perform in practice? To this end one would need to implement it and try it on a couple of test problems. Second, is the method also stable for $d>30$? We conjecture that this is the case, but for reasons of practical relevance we did not further investigate this.

Future work is related to implementing the method in an efficient way with an adaptive mesh-refinement strategy for each component $x_{k}$, $k=1,...,n_x$. Once the scheme is implemented, also the development of a NLP solver that is tailored to the solution of the unconstrained problem \eqref{eqn:NLP}, where $M\gg N$, coupled with a computational exploitation of the hierarchical discretization-matrices $\bP_{\dot{y},y,z}$ from former mesh-refinements in the linear algebra kernels, appears attractive for future research.

\FloatBarrier

\bibliography{OCP_disc_Bib}
\bibliographystyle{plain}

\end{document}